\documentclass{elsarticle}

\usepackage{amsmath}
\usepackage{bm}
\usepackage{tikz}
\usepackage{pgfplots}

\usepackage{hyperref}

\newtheorem{thm}{Theorem}
\newtheorem{lem}[thm]{Lemma}
\newproof{pf}{Proof}

\journal{arXiv.org} 

\graphicspath{{pictures/}}

\begin{document}

\begin{frontmatter}

\title{Two-level schemes for the advection equation}

\author[nsi,univ]{Petr N. Vabishchevich\corref{cor}}
\ead{vabishchevich@gmail.com}

\address[nsi]{Nuclear Safety Institute, Russian Academy of Sciences, 52, B. Tulskaya, Moscow, Russia}
\address[univ]{North-Eastern Federal University, 58, Belinskogo, Yakutsk, Russia}

\cortext[cor]{Corresponding author}

\begin{abstract}
The advection equation is the basis for mathematical models of continuum mechanics.
In the approximate solution of nonstationary problems it is necessary to inherit main properties of the conservatism and monotonicity of the solution.
In this paper, the advection equation is written in the symmetric form, where the advection operator is the half-sum of advection operators in conservative (divergent) and non-conservative (characteristic) forms.
The advection operator is skew-symmetric.
Standard finite element approximations in space are used.
The standart explicit two-level scheme for the advection equation is absolutly unstable. New conditionally stable regularized schemes are constructed,
on the basis of the general theory of stability (well-posedness) of operator-difference schemes, the stability conditions of the explicit Lax-Wendroff scheme are established.
Unconditionally stable and conservative schemes are implicit schemes of the second (Crank-Nicolson scheme) and fourth order. The conditionally stable implicit Lax-Wendroff scheme is constructed.
The accuracy of the investigated explicit and implicit two-level schemes for an approximate solution of the advection equation is illustrated by the numerical results of a model two-dimensional problem.
\end{abstract}

\begin{keyword}
Advection equation \sep two-level scheme \sep difference scheme stability \sep Lax-Wendroff scheme \sep Pade approximation \sep mass matrix diagonalization
\MSC[2010] 35L65 \sep 35L90 \sep 65M12 \sep 65M60 \sep 76M10
\end{keyword}


\end{frontmatter}

\section{Introduction}

Mathematical models of continuum mechanics \cite{Batchelor,LandauLifshic1986} describe the transport of scalar and vector quantities due to advection.
In particular, the basic equation of hydrodynamics is the continuity equation.
Advective transfer causes the fulfillment of conservation laws \cite{Godunov, LeVeque}.
In addition, there are properties of the positivity and monotonicity of the solution.
Such important properties of the differential problem must be inherited when passing to the discrete problem \cite{Anderson, Wesseling}.

For spatially approximation, conservative approximations are constructed on basis of using the conservative (divergent) form of the advection equation.
Most naturally such technology is implemented  when using the integro-interpolation method (balance method) on regular and irregular grids \cite{Samarskii1989}, in the control volume method \cite{LeVeque, Versteeg}.
The construction of monotonic approximations is discussed in many papers (see, for example, \cite{kulikovskii2000mathematical, HundsdorferVerwer2003, Kuzmin}).
In \cite{MortonKellogg1996, SamarskiiVabischevich1999a} standard linear approximations are considered for convection-diffusion problems.

Currently, the main computing technology for solving applied problems is the finite element method \cite{Guermond, Larson}.
It is widely used in computational fluid dynamics \cite{Donea, Zienkiewicz}.
Monotonization of the solution is achieved by using various linear and nonlinear variants of stabilization techniques.
It should be noted that the standard formulation of the equations of continuous medium mechanics in the conservative or non-conservative form is poorly suited for applying finite element approximations, for which the Hilbert spaces are natural.

Separate attention deserves the problems of constructing and investigating approximations in time.
When solving boundary value problems for partial differential equations, two-level schemes ($\theta$-method, schemes with weights) are traditionally widely used \cite{HundsdorferVerwer2003, Ascher2008, LeVeque2007}.
Research of this schemes can be based on the general theory of stability (well-posedness) of operator-difference schemes \cite{Samarskii1989, SamarskiiMatusVabischevich2002}.
In particular, unimprovable (coinciding necessary and sufficient) stability conditions can be used, which are formulated as operator inequalities in a finite-dimensional Hilbert space.
The achieved level of theoretical stability research allows us to abandon the heuristic methods widely used in computational fluid dynamics to study the stability of difference schemes: the von Neumann method for stability analysis, Fourier analysis, the principle of frozen coefficients, the consideration of the problem without taking into account the boundary conditions.

In this paper, the standard finite-element approximation in space is used for the nonstationary equation.
The advection equation is written in the so-called symmetric form \cite{SamarskiiVabischevich1999a}, when the advection operator is the half-sum of the advection operators in the conservative (divergent) and non-conservative (characteristic) forms.
Thus the continuum mechanics equations are written using the SD (Square root from Density) variables \cite{VabForm, form1}.
In this case the corresponding conservation laws are a direct consequence of the skew-symmetry of the advection operator.
The conservativeness property is related to the preservation of the norm of the solution of the non-stationary advection equation, with stability with respect to the initial data.
This property takes place not only for the solution, but also for some of its transformations.
In this case, we are talking about the property of multiconservativeness.
The principal point is related with the fact that the most important skew-symmetry property of the advection operator is inherited for finite element approximations.

For the advection equation an explicit two-level scheme as well as all schemes with a weight lower than 0.5, are aboslutly unstable.
At the same time, extensive computational practice aims us to use the explicit schemes in these problems.
In such schemes, the conditional stability (the Courant-Friedrichs-Lewy (CFL) condition) with a time-step limited by the Courant number is provided in fact by refusing the skew-symmetry of the advection operator --- the use of dissipative approximations.
In addition, the standard explicit schemes are not conservative.

We construct conditionally stable schemes for the advection equation based on the principle of regularization of operator-difference schemes \cite{Samarskii1967, Vabishchevich2014}, when stability is provided by a small perturbation.
The explicit second-order Runge-Kutta scheme \cite {Butcher2008, HairerWanner2010} is considered.
In this context, the classical explicit Lax-Wendroff scheme \cite{lax1964difference} is also considered as regularization of the second-order Runge-Kutta scheme --- the perturbation of the advection operator squared.
Using the stability criteria of two-level operator-difference schemes, the stability conditions of regularized schemes and the explicit Lax-Wendroff scheme are obtained.
Effective computational implementation of explicit schemes using finite element approximation in space is provided by using the diagonalization procedures of the mass matrix (mass-lumping procedure) \cite{Thomee2006, Cohen}.

Of greatest interest are implicit two-level schemes for the advection equation, which belong to the unconditionally stable class.
The classical Crank-Nicolson scheme has a second-order of accuracy, is unconditionally stable and multiconservative.
For the advection problems under the consideration, we can use a scheme of the fourth-order of accuracy, which is also unconditionally stable and multiconservative.
A certain drawback of this scheme is associated with the need to use the lumping procedure.
An implicit version of the Lax-Wendroff scheme is proposed, which is conditionally stable, but has a higher accuracy than the explicit Lax-Wendroff scheme and does not require the diagonalization of the mass matrix.

The paper is organized as follows. A model two-dimensional problem for the advection equation is formulated in Section 2. 
Approximation in space is constructed using Lagrangian finite elements, the main properties of the problem solution are noted.
In Section 3, we consider known and new explicit difference schemes for the advection equations, and investigate the stability conditions.
Central for this work is Section 4.
Implicit schemes, their stability and conservatism are studied here.
In Section 5, numerical experiments on the accuracy of explicit and implicit schemes are discussed for the model IBV problem. The results of the work are summarized in Section 6.

\section{Problem statement}

In a bounded two-dimensional domain $\Omega$, we consider the advection equation written in the symmetric form.
A standard finite element approximation in space is used.
The problem of constructing approximations in time is formulated in such a way that the approximate solution inherits the basic properties of the solution of the differential problem.

\subsection{Differential problem}
The Cauchy problem is considered in the domain $\Omega$ ($\bm x = (x_1, x_2) \in \Omega$)
\begin{equation}\label{1}
 \frac{d w}{d t} + \mathcal{A} w = 0,
 \quad 0 < t \leq T, 
\end{equation} 
\begin{equation}\label{2}
 w(0)= w^0,
\end{equation}
using notation $w(t) = w(\bm x,t)$.
The operator of advection (convective transport) $\mathcal{A}$ is assumed to be constant and is written in the symmetric form:
\begin{equation}\label{3}
 \mathcal{A} w =  \frac{1}{2} {\rm div}  (\bm v \, w) + 
 \frac{1}{2} \bm v \cdot {\rm grad}   \, w.
\end{equation}
Thus, we take the half-sum of the transfer operator in the divergent (conservative part with ${\rm div} (\bm v \, w)$) and non-divergent (characteristic part with $ \bm v \cdot {\rm grad} \, w $) forms \cite{SamarskiiVabischevich1999a, zang1991rotation}.
The convective transport is determined by the velocity of the medium $\bm v (\bm x)$, and the no-permeability condition on the boundary of the domain
\begin{equation}\label{4}
 (\bm v \cdot \bm \nu ) = 0, \quad \bm x \in \partial \Omega ,  
\end{equation}
where $\bm \nu$ is the outer normal to the  boundary $\partial \Omega$. 

The standard continuity equation has a divergent form:
\begin{equation}\label{5}
  \frac{\partial \varrho }{\partial t} +
  {\rm div} (\varrho {\bm v} ) = 0 , 
\end{equation} 
where $\varrho$ is the density.  
From (\ref{4}), (\ref{5}), by direct integration, we obtain the law of mass conservation:
\[
 m (t) = \mathrm{const} ,
 \quad m = \int_{\Omega}  \varrho(\bm x, t) d \bm x .
\]
When orienting to finite elemental approximations in space, it is convenient to use of SD (Square root from Density) variables for the equations of hydrodynamics \cite{VabForm}.
Let $w = \varrho^{1/2}$, then equation (\ref{5}) is written in the form (\ref{1}), (\ref{3}).

In the Hilbert space $H = L_2(\Omega)$, we define the scalar product and norm in the standard way:
\[
  (w,u) = \int_{\Omega} w({\bm x}) u({\bm x}) d{\bm x},
  \quad \|w\| = (w,w)^{1/2} .
\] 
Under constraint (\ref{4}), the convective transfer operator is skew-symmetric:
\begin{equation}\label{6}
 \mathcal{A}  = -  \mathcal{A}^* ,
\end{equation} 
and therefore, in particular,
\[
 (\mathcal{A} u, u) = 0 .
\] 

Using  the scalar product of equation (\ref{1}) and $w$, taking into account (\ref{2}), (\ref{6}), we obtain
\begin{equation}\label{7}
 \|w(t)\| = \|w^0\| ,
 \quad 0 < t \leq T .
\end{equation} 
This relation with respect to the continuity equation (\ref{5}) expresses the law of mass conservation.
The property (\ref{7}) for the solution of the the Cauchy problem (\ref{1}), (\ref{2}) is associated with the non-dissipativity of the system.
The equation (\ref{1}) for (\ref{6}) has many other conservation laws.
If the constant operator $\mathcal{B}$ is permutable with $\mathcal{A}$, then we have the equality
\begin{equation}\label{8}
 \|\mathcal{B} w(t)\| = \|\mathcal{B} w^0\| ,
 \quad 0 < t \leq T .
\end{equation} 
In particular, equality (\ref{8}) holds for $\mathcal{B} = f(\mathcal{A})$.

\subsection{Finite element approximation in space} 

To solve numerically the problem (\ref{1})--(\ref{3}), we employ finite element approximations in space (see, e.g., \cite{Thomee2006,brenner2008mathematical}). 
For (\ref{3}), we define the bilinear form
\[
 a(w,u) = \int_{\Omega } \left ( \frac{1}{2} {\rm div}  (\bm v \, w) \, u + 
 \frac{1}{2} \bm v \cdot {\rm grad} \, w \, u \right )  d {\bm x} .
\] 
By (\ref{4}), we have
\[
 a(w,u) = - a(u,w),
 \quad a(w,w) = 0 .   
\]
Define the subspace of finite elements $V^h \subset H^1(\Omega)$ and the discrete operator $A$ as
\[
(A w, u) = a(w,u),
\quad \forall \ w, u \in V^h . 
\]
The operator $A$ acts on the finite dimensional space $V^h$ 
and, similarly to (\ref{8}), is skew-symmetric: 
\begin{equation}\label{9}
  A =-  A^* . 
\end{equation} 
In finite element approximation, the main property of the advection operator, namely, its  skew-symmetry, is inherited and, as a consequence, there is energy neutrality (the equality $(Aw,w) = 0$).

The Cauchy problem (\ref{1}), (\ref{2}) is associated with the problem
\begin{equation}\label{10}
 \frac{d u}{d t} + A u = 0,
 \quad 0 < t \leq T, 
\end{equation} 
\begin{equation}\label{11}
 u(0)= u^0 ,
\end{equation} 
for $u(t) \in V^h$, where $u^0 = P w^0$ with $P$ denoting $L_2$-projection onto $V^h$.
Similarly (\ref{7}), for the solution of the problem (\ref{10}), (\ref{11}), the conservative property is established:
\begin{equation}\label{12}
 \|u(t)\| = \|u^0\| ,
 \quad 0 < t \leq T .
\end{equation} 
The multiconservative property (see (\ref{8})) is associated with the equality
\begin{equation}\label{13}
 \|B u(t)\| = \|B u^0\| ,
 \quad 0 < t \leq T ,
\end{equation} 
provided that
\[
 \frac{d}{d t} B = B \frac{d}{d t} ,
 \quad B A = A B . 
\] 
When the problem (\ref{10}), (\ref{11}) is approximated in time, it is necessary to focus on the fulfillment of the properties (\ref{12}) and (\ref{13}) at separate time levels.

\subsection{Approximation in time} 

Let, for simplicity, $\tau$ be a step of a uniform grid in time such that $y^n = y(t^n), \ t^n = n \tau$, $n = 0,1, ..., N, \ N\tau = T$.
To solve numerically the problem (\ref{10}), (\ref{11}), we use explicit and implicit two-level schemes, when the solution at the new level $y^{n+1} $ is determined by the previously found solution $y^{n}$, $n = 0, 1, ..., N-1$.

When approximating in time, we focus, first of all, on the scheme stability.
With respect to our problem, stability will be ensured by the following estimate of the solution at each time level:
\begin{equation}\label{14}
 \|y^{n+1}\| \leq \exp(\mu \tau)  \|y^{n}\|,
 \quad \mu =  \mathrm{const} ,
 \quad n = 0,1, ..., N-1 .  
\end{equation} 

The most favorable situation is associated with the fulfillment of the equality
\begin{equation}\label{15}
 \|y^{n+1}\| = \|y^{n}\|,
 \quad n = 0,1, ..., N-1 .   
\end{equation} 
In this case, not only stability is ensured, but also the conservatism of the solution takes place (see (\ref{12})). Multiconservatism (see (\ref {13})) is due to the fulfillment of the equality
\begin{equation}\label{16}
 \|B y^{n+1}\| = \|B y^{n}\|,
 \quad n = 0,1, ..., N-1 ,  
\end{equation} 
for some operators $B$.

Conservative approximations are related to the property of time reversibility for the problem (\ref{9})--(\ref{11}).
The computational algorithm is reversible in time if we calculate the solution at the level $t^{n+1}$ and then change the transfer velocity ($\bm v (\bm x) $) to the opposite (by $-\bm v (\bm x )$), then at the next step in time we get the solution that coincides with the solution at the level $t^{n}$.

\section{Explicit schemes}

The explicit scheme for the problem (\ref{9})--(\ref{11}) has the first-order error in time and is absolutely unstable.
A conditionally stable scheme can be constructed on the basis of the regularization principle of difference schemes.
We also consider a scheme of second-order approximation, namely, the Lax-Wendroff scheme.

\subsection{Necessary and sufficient conditions for the stability of schemes with weights}

Strict results on the stability of difference schemes in finite-dimensional Hilbert spaces are obtained in the theory of stability (well-posedness) of operator-difference schemes \cite{Samarskii1989, SamarskiiMatusVabischevich2002}.
With respect to the subject of our investigation, we give the best possible conditions for the stability of a standard two-level scheme with weights ($\theta$-method) to be unimprovable (matching necessary and sufficient) for the solution of the problem (\ref{10}), (\ref{11}).

When passing from the level $t^n$ to the level $t^{n+1}$, the approximate solution of the problem (\ref{10}), (\ref{11}) is determined from the equation
\begin{equation}\label{17}
 \frac{y^{n+1} - y^n}{\tau} + C (\theta y^{n+1} + (1 - \theta) y^n ) = 0,
 \quad n = 0,1,..., N-1 .  
\end{equation} 
with some constant operator $C$.
In the simplest case, $C=A$.
The initial condition is
\begin{equation}\label{18}
 y^0 = u^0.
\end{equation}
If $\theta=0$, we have then explicit scheme, for $\theta=1$, we obtain the fully implicit scheme and for $\theta=0.5$, the Crank-Nicolson scheme appears.
The stability criterion is formulated as follows.
\begin{lem}
Condition
\begin{equation}\label{19}
 (C y,y) + \left (\theta - \frac{1}{2} \right )  \tau \|C y\|^2 \geq 0  
\end{equation} 
is necessary and sufficient for the stability of the scheme (\ref{17}), (\ref{18}) in $H$, and for the solution  the level-wise inequality (\ref{14}) holds with $\mu=0$.
\end{lem}

The proof of this statement is given, for example, in the book \cite{SamarskiiGulin1973}.
The condition (\ref{19}) can be written in the form of the operator inequality
\[
 C + \left (\theta - \frac{1}{2} \right ) \tau C^* C \geq 0.  
\]

\subsection{Explicit schemes of the first and second order of approximation} 

It is natural to start with the simplest explicit scheme for the problem (\ref{10}), (\ref{11}).
In this case, $C = A, \ \theta = 0$ in (\ref {17}), i.e.
\begin{equation}\label{20}
 \frac{y^{n+1} - y^n}{\tau} + A y^n = 0,
 \quad n = 0,1,..., N-1 .  
\end{equation}  

For the explicit scheme, the stability criterion (\ref{19}) takes the form
\begin{equation}\label{21}
 (C y,y) - \frac{\tau}{2}  \|C y\|^2 \geq 0 .
\end{equation} 
In the case of a skew-symmetric operator $C$,  inequality (\ref{21}) is not satisfied for any $\tau>0$.
Therefore, the explicit scheme (\ref{18}), (\ref{20}) is absolutely unstable.

Among two-level schemes, explicit schemes of the second-order approximation are deserved separate consideration.
Instead of (\ref{20}) we will use the explicit second-order Runge-Kutta scheme:
\begin{equation}\label{22}
 \frac{y^{n+1} - y^n}{\tau} + A y^n - \frac{\tau}{2} A^2 y^n = 0,
 \quad n = 0,1,..., N-1 .  
\end{equation}  
In this case, we have
\[
 C = A - \frac{\tau}{2} A^2 .
\] 
Taking into account (\ref{9}), we get
\[
\begin{split}
 C^* C & = \left (A^* + \frac{\tau}{2} A^* A \right ) \left (A + \frac{\tau}{2} A^* A \right ) \\
& = \frac{\tau^2}{4} (A^* A)^2 + A^* A.
\end{split}
\] 
In view of this, the stability condition (\ref{21}) takes the form
\[
 - \frac{\tau^3}{8} \|A^2 y\|^2 \geq 0 .
\]
Again, we cannot specify $\tau > 0 $ such that the stability of the scheme (\ref{22}) holds.

Thus, we can formulate the following result.
\begin{thm}
The explicit scheme of the first-order approximation (\ref{18}), (\ref{20}) and the explicit scheme of the second-order approximation (\ref{18}), (\ref{22}) are absolutly unstable in $H$.
\end{thm}

\subsection{Conditionally stable non-standard scheme} 

We note the possibility of constructing a conditionally stable scheme based on the correction of the explicit scheme (\ref{18}), (\ref{20}).
We consider more weak stability requirements, allowing the growth of the solution norm in accordance with (\ref{14}) with the choice of some positive constant $\mu > 0$.

The construction of the considering approximation in time is based (see, for example, \cite{afanas2013unconditionally, VabExact}) on the solution representation of the problem (\ref{10}), (\ref{11}) in the form
\begin{equation}\label{23}
 u(t) = \exp(\mu t) v(t) .
\end{equation} 
For $v(t)$, from (\ref{10}), (\ref{11}), (\ref{23}), we have the Cauchy problem
\begin{equation}\label{24}
 \frac{d v}{d t} + (A + \mu I) v = 0,
\end{equation} 
\begin{equation}\label{25}
 v(0) = u^0 ,
\end{equation} 
where $I$ is the identity operator.

For the approximate solution of the problem (\ref{24}), (\ref{25}), the explicit difference scheme is used
\begin{equation}\label{26}
 \frac{v^{n+1} - v^n}{\tau} + (A + \mu I) v^n = 0,
 \quad n = 0,1,..., N-1 ,  
\end{equation}  
\begin{equation}\label{27}
 v^0 = u^0 .
\end{equation}
Taking into account relation
\[
 y^{n} = \exp(\mu t^n) v^n,
 \quad n = 0,1, ..., N 
\]  
scheme (\ref{26}), (\ref{27}) corresponds to the use of
\begin{equation}\label{28}
 \frac{\exp(-\mu \tau) u^{n+1} - u^n}{\tau} + (A + \mu I) u^n = 0,
 \quad n = 0,1,..., N-1 , 
\end{equation}
with the initial condition (\ref{18}).
Such schemes belong to the class of non-standard \cite{mickens1994nonstandard,mickens2002nonstandard}.
The level-wise estimate
\[
 \|v^{n+1}\| \leq \|v^n\|
\] 
for the solution of the difference scheme (\ref {26}), (\ref{27}) corresponds to the estimate (\ref{14}) for the scheme (\ref{18}).

\begin{thm}
The explicit scheme (\ref{18}), (\ref{28}) is stable in $H$ for
\begin{equation}\label{29}
 \tau \leq \frac{2 \mu }{\mu^2 + \|A\|^2} ,
 \quad \mu > 0 . 
\end{equation} 
In these conditions, for the problem solution the estimate (\ref{14}) is satisfied.
\end{thm}

\begin{pf}
It suffices to formulate the stability conditions of the scheme
(\ref{26}), (\ref{27}).
In this case $C = A + \mu I$ and the stability criterion (\ref{21}) gives
\[
 \mu I - \frac{\tau}{2} (\mu^2 I + A^* A) \geq 0 .
\] 
This inequality will be satisfied with time-step constraints
This inequality will be satisfied with constraints of time step (\ref{29}).
\end{pf}

For the advection problems under the consideration, an acceptable time step is associated with the Courant condition, when
\[
 \tau \leq \mathcal{O} (\|A\|^{-1}) .
\] 
The time step limitation (\ref{29}) is much more stiff:
$\tau \leq \mathcal{O} (\|A\|^{-2})$.
Therefore, the explicit scheme (\ref{18}), (\ref{28}) is not suitable for computational practice.

\subsection{Regularized schemes}

The regularization principle for difference schemes \cite{Samarskii1989,Vabishchevich2014} provides great opportunities in constructing difference schemes of a prescribed quality. 
The standard approach to the construction of stable schemes on the basis of the regularization principle is associated with the introduction of additional terms (regularizers) in operators of an origional (generating) difference scheme.

Instead of the explicit scheme (\ref{20}), we will use the regularized scheme
\begin{equation}\label{30}
 \frac{y^{n+1} - y^n}{\tau} + A y^n + \frac{\tau}{2} D y^n = 0,
 \quad n = 0,1,..., N-1 .  
\end{equation}  
with the regularizer $D = D^* > 0$.
Let us formulate the stability conditions for this scheme.

For the scheme (\ref{30}), we have
\[
 C = A + \frac{\tau}{2} D .
\] 
Taking into account the skew-symmetry of the operator $A$, we obtain
\[
\begin{split}
 C^* C & = \left (A^* + \frac{\tau}{2} D \right ) \left (A + \frac{\tau}{2} D \right ) \\
& = \frac{\tau^2}{4} D^2 + A^* A.
\end{split}
\] 
The stability condition (\ref{21}) gives
\begin{equation}\label{31}
 \frac{\tau^2}{4} D^2 \leq  D - A^* A . 
\end{equation}
Therefore, we can rely on the conditional stability of the scheme (\ref{30}) for $D \geq A^* A$.

Let $\lambda_{\max}$ be the maximum eigenvalue of the spectral problem
\begin{equation}\label{32}
 D^2 = \lambda (D - A^* A ) ,
\end{equation} 
and $\lambda_{\max} > 0$. Then the inequality (\ref{31}) will be satisfied for
\begin{equation}\label{33}
 \tau \leq \tau_0,
 \quad \tau_0 = \frac{2}{\lambda_{\max}^{1/2}}.
\end{equation} 
Thus, we can formulate the stability conditions.

\begin{thm}\label{t-4}
The explicit scheme (\ref{18}), (\ref{30}) for $D = D^* \geq A^* A$ is stable in $H$ if the time step satisfies the constraints (\ref{33}), where $\lambda_{\max}$ is the maximum eigenvalue of the spectral problem (\ref{32}).
\end{thm}

We highlight some interesting possibilities for the choice of the regularizer $D$.
The simplest version is related to the regularizer
\begin{equation}\label{34}
 D = \beta A^* A,
 \quad \beta > 1. 
\end{equation} 
The conditions (\ref{32}), (\ref{33}) give the following restrictions on the time step:
\begin{equation}\label{35}
 \tau \leq 2 \frac{(\beta-1)^{1/2}}{\beta} \frac{1}{\|A\|} .
\end{equation} 
Thus, the explicit scheme (\ref{18}), (\ref{30}), (\ref{34}) is conditionally stable under the time step constraints of Courant type (see (\ref{35})).
In the special case, for $\beta = 2$, we have
\[
 \tau \leq \frac{1}{\|A\|} . 
\] 

The regularized scheme (\ref{18}), (\ref{30}), (\ref{34}) has the first-order approximation in time.
It is natural to expect a higher accuracy when choosing $D \approx A^* A $, i.e. when the scheme (\ref{30}) is close to the explicit second-order Runge-Kutta scheme (\ref{22}).
Here we can indicate two variants of such schemes.

The first variant is related with choosing the regularizer $D$ in the form
\begin{equation}\label{36}
 D = (1+\beta \tau) A^* A,
 \quad \beta > 0. 
\end{equation} 
The stability condition (\ref{21}) for this case takes the form
\[
 \frac{\tau^2}{4} (1 + \beta \tau)^2 (A^* A)^2 \leq \beta \tau A^* A.
\] 
It will be satisfied for
\begin{equation}\label{37}
 \tau \leq 4 \beta \frac{1}{\|A\|^2} .
\end{equation} 
The restrictions (\ref{37}) for the scheme (\ref{18}), (\ref{30}), (\ref{36}) are substantially more strong than the restrictions (\ref{35}) for the scheme (\ref{18}), (\ref{30}), (\ref{34}).
This disadvantage is compensated by the fact that we have the second-order approximation in time instead of the first order.

Let us now consider the second variant of the regularized scheme of the second-order approximation.
The most well-known modification of the absolutely unstable explicit second-order Runge-Kutta scheme (\ref{18}), (\ref{22}) is the Lax-Wendroff scheme \cite{lax1964difference, leveque2002finite}.
This scheme is actually based on replacing the operator $-A^2 = A^* A$ in the explicit second-order Runge-Kutta scheme by a self-adjoint non-negative operator $Q$ close to it.

We define the new bilinear form
\[
 q(w,u) = \int_{\Omega } \left ( \frac{1}{2} {\rm div}  (\bm v \, w)  + 
 \frac{1}{2} \bm v \cdot {\rm grad} \, w \right )  
\left ( \frac{1}{2} {\rm div}  (\bm v \, u) + 
 \frac{1}{2} \bm v \cdot {\rm grad} \, u \right ) 
d {\bm x} .
\] 
By virtue of this definition, we have
\[
 q(w,u) = q(u,w),
 \quad q(w,w) \geq  0 .   
\]
Define the discrete operator $Q$ as
\begin{equation}\label{38}
(Q w, u) = q(w,u),
\quad \forall \ w, u \in V^h .  
\end{equation} 
The regularized scheme (\ref{18}), (\ref{30}) for $D = Q$ is a finite element version of the explicit Lax-Wendroff scheme.
The stability conditions for $ D \geq A ^ * A $ on a finite element space $ V ^ h $ are formulated in the theorem~\ref{t-4}. 

\subsection{Computational implementation of explicit schemes}

The explicit schemes under the consideration are explicit only in form, according to the time approximation.
The computational implementation of explicit finite-element approximations, unfortunately, is connected with the solution of systems of linear equations on a new level in time.
The standard approach is related to the correction of approximations based on mass lumping (see, for example, \cite{Thomee2006}).

The use of certain schemes for the Cauchy problem (\ref{10}), (\ref{11}) is based on the solution of matrix problems for finding an approximate solution at a new time level.
We associate with the differential-operator equation the corresponding system of ordinary differential equations.

We consider a standard quasi-uniform triangulation of the domain $\Omega$ into triangles (or tetrahedra in 3-D). 
Let $ x_i, \ i = 1,2, ..., N_h$ be vertices of this triangulation. We introduce the finite dimensional space $V^h \subset H^1(\Omega)$ of continuous functions that are liner over each finite element, see, e.g. \cite{Thomee2006}. 
As a nodal basis we take the standard \emph{hat} function $\chi_i(x) \in V^h, \ i = 1,2, ..., N_h$.
Then for $v \in V^h$, we have the representation
\[
 v(x) = \sum_{i=i}^{N_h} v_i \chi_i( x).
\] 
where $v_i = v(x_i), \ i = 1,2, ..., N_h$.

From (\ref{10}), we obtain the equation
\begin{equation}\label{39}
 M \frac{d z}{d t} + K z = 0,
\end{equation}
where $z(t)$ is the vector of unknowns $z_i(t), \ i = 1,2, ..., N_h$.
Here $M = (m_{ij})$ is the mass matrix and $K = (k_{ij})$ is the stiffness matrix, and
\[
 m_{ij} = (\chi_i, \chi_j), 
 \quad k_{ij} = a(\chi_i, \chi_j), 
 \quad i,j = 1,2, ..., N_h . 
\]
If we set $u = M^{1/2} z$, then equation (\ref{39}) is written in the form (\ref{9}), (\ref{10}) for
\begin{equation}\label{40}
 A = M^{-1/2} K M^{-1/2} ,
 \quad M = M^* > 0, \quad K = - K^* . 
\end{equation} 
Thus, we can construct approximations in time for the equation (\ref{39}) on the basis of approximations for equation (\ref{10}), taking into account that $u = M^{1/2} z$ and (\ref {40}).

For example, the explicit scheme (\ref {20}) corresponds to the scheme
\begin{equation}\label{41}
 M \frac{z^{n+1} - z^n}{\tau} + K z^n = 0,
 \quad n = 0,1,..., N-1 .  
\end{equation} 
It is necessary to solve the problem
\[
  M z^{n+1} = M z^n - \tau K z^n
\] 
at a new level in time.
Because of this, the scheme (\ref{41}) is explicit in terms of time approximation, but implicit in terms of the computational implementation, since the symmetric matrix $M$ is off-diagonal.

In order to provide an explicit computational implementation, various diagonalization procedures for the mass matrix are used:
\[
 M \longrightarrow \widetilde{M} .
\] 
In the simplest mass lumping procedure \cite{Thomee2006} we have
\[
 \widetilde{M} =  \mathrm{diag} \{ \widetilde{m}_1, \widetilde{m}_2, ..., \widetilde{m}_{N_h} \} ,
 \quad \widetilde{m}_{i} = \sum_{j=1}^{N_h} m_{ij} ,
 \quad i = 1,2, ... , N_h . 
\] 
Instead of (\ref {39}), we seek an approximate solution of the Cauchy problem for the equation
\[
 \widetilde{M} \frac{d z}{d t} + K z = 0 .
\]
This equation corresponds to setting
\begin{equation}\label{42}
 A = \widetilde{M}^{-1/2} K \widetilde{M}^{-1/2} ,
 \quad \widetilde{M} = \widetilde{M}^* > 0, \quad K = - K^* , 
\end{equation} 
in the problem (\ref{9})--(\ref{11}).

\section{Implicit schemes} 

Unconditionally stable schemes for the advection equation are built on the basis of implicit approximations.
It seems reasonable to employ the Crank-Nicolson scheme, which has the second order of accuracy in time and has conservative properties.
In addition, schemes of higher accuracy are outlined, which have the fourth-order approximation.
The implicit version of the Lax-Wendroff scheme is highlighted.

\subsection{The Crank-Nicolson scheme}

Among the implicit schemes for the problem (\ref{10}), (\ref{11}), the most important is the symmetric scheme whith $C = A, \ \theta = 0.5$ in (\ref{17}).
In particular, it has the best SM (Spectral Mimetic) properties in the class of schemes with the skew-symmetric operator \cite{Vabischevich2011}.
In this case, the solution is determined from the equation
\begin{equation}\label{43}
 \frac{y^{n+1} - y^n}{\tau} + A \frac{y^{n+1} + y^n}{2} = 0,
 \quad n = 0,1,..., N-1 .  
\end{equation}
By virtue of (\ref{19}), the scheme (\ref{18}), (\ref{43}) is absolutely stable and approximates the problems (\ref{10}), (\ref{11}) with the second order in $\tau$.

The Crank-Nicolson scheme is nondissipative. To show this fact, it is sufficient to multiply the scalar equation (\ref{43}) by $\tau (y^{n+1} + y^n)$.
Taking into account the skew-symmetry of the operator $A$, we have
\[
 \|y^{n+1}\|^2 - \|y^{n}\|^2 = 0,
 \quad n = 0,1, ..., N-1 . 
\]
Taking into account (\ref{18}), we arrive at (\ref{15}).
Also there is the multiconservative property, see (\ref{13}).
For any constant operator $B$, which commutes with $A$, (\ref{16}) is satisfied.
We immediately see that the scheme (\ref{18}), (\ref{43}) generates a computational algorithm that is time reversible.
The result of our consideration is the following theorem.

\begin{thm}\label{t-5}
The scheme (\ref{18}), (\ref{43}) provides a time-reversible computational algorithm and is unconditionally stable in $H$.
The scheme is conservative and multiconservative in the sense of the fulfillment of the equalities (\ref{15}) and (\ref{16}) for all operators $B$ that are constant and permutable with $A$.
\end{thm}

\subsection{Scheme of higher accuracy order}

The Crank-Nicolson scheme (\ref{43}) corresponds to the use of the following Pade approximant for the exponential:
\[
 \exp(-z) = \frac{1 - \frac{1}{2} z }{1 + \frac{1}{2} z} + \mathcal{O}(z^3) . 
\] 
When considering more accurate approximations, we can consider
\[
 \exp(-z) = \frac{1 - \frac{1}{2} z + \frac{1}{12} z^2}{1 + \frac{1}{2} z + \frac{1}{12} z^2} + \mathcal{O}(z^5) . 
\] 
The corresponding difference scheme ($z = \tau A$) has the form
\begin{equation}\label{44}
 \left (I + \frac{1}{12} \tau^2 A^2 \right ) \frac{y^{n+1} - y^n}{\tau} + A \frac{y^{n+1} + y^n}{2} = 0,
 \quad n = 0,1,..., N-1 .  
\end{equation}  

We set
\[
 E = I + \frac{1}{12} \tau^2 A^2 ,
\] 
and apply scalar multiplication of equation (\ref{44}) by $\tau(y^{n+1} + y^n)$. This leads to the equality
\[
 (Ey^{n+1}, y^{n+1}) =  (Ey^{n}, y^{n}) ,
 \quad n = 0,1,..., N-1 .   
\]
For $E > 0$, this equality, taking into account $EA = AE$, implies the stability estimate, the conservativeness and multiconservative properties.
The positivity of the operator $E$ is satisfied for at not very large steps in time.
Really,
\[
 E = I - \frac{1}{12} \tau^2 A^*A \geq  I - \frac{1}{12} \tau^2 \|A\|^2 I . 
\]
Thus $E > 0$ for
\begin{equation}\label{45}
 \tau < \frac{2 \cdot 3^{1/2}}{\|A\|} . 
\end{equation} 

The restrictions (\ref{45}) can be removed by modifying the proof.
We write the scheme (\ref{44}) in the form
\begin{equation}\label{46}
 S_1 y^{n+1} = S_2 y^{n} ,
 \quad n = 0,1,..., N-1 , 
\end{equation} 
with notation
\begin{equation}\label{47}
 S_1 = I + P + \frac{1}{3} P^2,
 \quad S_2 = I - P + \frac{1}{3} P^2,
 \quad P = - P^*,  
\end{equation} 
and ${\displaystyle P = \frac{\tau }{2} A}$.
The first question is related with the proof of the existence $S_1^{- 1}$.
Necessary and sufficient condition for the existence of the inverse operator $S_1^{-1}$ is (see, for example, \citep{Anderson, Lusternik}) the fulfillment of the inequality
\[
 \|S_1 y\| \geq \delta \|y\|,
 \quad \delta > 0 , 
\] 
wherein $\|S_1^{-1}\| \leq \delta^{-1}$.

\begin{lem}\label{t-6}
For the operator $S_1$, which is defined according to (\ref{47}), the inequality is satisfied
\begin{equation}\label{48}
 \|S_1 y\| \geq  \|y\| . 
\end{equation} 
\end{lem}

\begin{pf}
The inequality (\ref{48}) is equivalent to the inequality
\[
 (S_1^* S_1 y, y) \geq (y,y) .
\] 
We have
\[
\begin{split}
 S_1^* S_1 & =  \left (I + \frac{1}{3} P^2 - P \right )  \left (I + \frac{1}{3} P^2 + P \right ) \\
 & = \left (I + \frac{1}{3} P^2 \right )^2 - P^2 
 = I - \frac{1}{3} P^2 + \frac{1}{9} P^4 \\
 & = I + \frac{1}{3} P^* P + \frac{1}{9} P^* P  P^* P \geq I 
\end{split}
\] 
taking into account the skew-symmetry of the operator $P$.
\end{pf}

Taking into account lemma \ref{t-6}, we can write (\ref{46}) as follows
\begin{equation}\label{49}
 y^{n+1} = S y^{n} ,
 \quad S =  S_1^{-1} S_2, 
 \quad n = 0,1,..., N-1 .
\end{equation} 
For the transition operator from one level in time to another, the following statement holds.

\begin{lem}\label{t-7}
The operator $S = S_1^{-l} S_2$ under the conditions (\ref{47}) is unitary, i.e.
\begin{equation}\label{50}
 S^* S = I .
\end{equation} 
\end{lem}

\begin{pf}
We have
\[
 S = \left (I + \frac{1}{3} P^2 + P \right )^{-1} \left (I + \frac{1}{3} P^2 - P \right ) , 
\] 
\[
 S^* = \left (I + \frac{1}{3} P^2 + P \right ) \left (I + \frac{1}{3} P^2 - P \right )^{-1} .
\] 
Taking into account the permutability of the operator factors on the right-hand side, we get the equality (\ref{50}).
\end{pf}

The noted unitarity property (\ref{50}) of the operator $S$ allows (\ref{49}) to come to the conservativity property (\ref{15}).
Similarly, the multiconservative property is established.
The result of our consideration is the following statement.

\begin{thm}\label{t-8}
The fourth-order accuracy scheme (\ref{18}), (\ref{44}) provides a time-reversible computational algorithm and is unconditionally stable in $H$.
The scheme (\ref{18}), (\ref{44}) is conservative and multiconservative in the sense of the equality (\ref{15}) and (\ref{16}) for all operators $B$ that are constant and permutable with $A$.
\end{thm}

In the computational implementation, we solve the Cauchy problem for equation (\ref{39}).
Because of this, the representation (\ref{40}) is used in equation (\ref{10}).
The scheme (\ref{44}) corresponds to the scheme
\begin{equation}\label{51}
\begin{split}
 \left (M + \frac{1}{12} \tau^2 K M^{-1} K \right ) \frac{z^{n+1} - z^n}{\tau} + & K \frac{z^{n+1} + z^n}{2} = 0, \\
 & n = 0,1,..., N-1 . 
\end{split}
\end{equation}  
Thus, for the transition to a new  time level, it is necessary to solve equation $Rz^{n+1} = r^n$ with the operator
\[
 R = M + \frac{1}{2} \tau K + \frac{1}{12} \tau^2 K M^{-1} K .
\] 
This makes the scheme (\ref{18}), (\ref{44}) practically useless.

As in the case of explicit schemes, the diagonalizing procedure of the mass matrix saves the situation.
In this case, instead of (\ref{40}), we use the representation (\ref{42}), and the scheme (\ref{45}) is replaced by the scheme
\begin{equation}\label{52}
\begin{split}
 \left (\widetilde{M} + \frac{1}{12} \tau^2 K \widetilde{M}^{-1} K \right ) \frac{z^{n+1} - z^n}{\tau} + & K \frac{z^{n+1} + z^n}{2} = 0, \\
 & n = 0,1,..., N-1 . 
\end{split}
\end{equation} 
The scheme (\ref{52}) is stable under constraint (\ref{45}) taking into account the fact that the representation (\ref{42}) holds.

\subsection{Implicit Lax-Wendroff scheme} 

On the base of the implicit scheme (\ref{44}), we can construct an implicit version of the Lax-Wendroff scheme.
Taking into account the notation introduced above, instead of (\ref{44}), we will use the scheme
\begin{equation}\label{53}
 \left (I - \frac{1}{12} \tau^2 Q \right ) \frac{y^{n+1} - y^n}{\tau} + A \frac{y^{n+1} + y^n}{2} = 0,
 \quad n = 0,1,..., N-1 .  
\end{equation} 
As in the case of the explicit scheme, the construction is based on replacing the operator $-A^2 = A^* A$ by the operator $Q$, which is defined by (\ref{38}).

We rewrite the scheme (\ref{53}) in the form
\[
 \left (I + \frac{1}{12} \tau^2 A^2 \right ) \frac{y^{n+1} - y^n}{\tau} + A \frac{y^{n+1} + y^n}{2} -  \frac{1}{12} \tau^2 R \frac{y^{n+1} - y^n}{\tau} = 0,
\] 
where
\[
 R = Q - A^* A .
\] 
The operator $R$ explicitly highlights the approximation error in space due to the replacement of $-A^2 $ by $Q$.

For a positive constant selfadjoint operator $E$, we define the Hilbert space $H_E$, where the scalar product and norm are defined as follows
\[
 (u,v)_E = (E u,v),
 \quad \|u\|_E = (u,u)_E^{1/2} .
\] 
Analogously to theorem \ref{t-8}, the following result is formulated.

\begin{thm}\label{t-9}
The implicit Lax-Wendroff scheme (\ref{18}), (\ref{53}) provides a time-reversible computational algorithm and is stable for the constraint
\begin{equation}\label{54}
 \tau < \frac{2 \cdot 3^{1/2}}{\|Q\|^{1/2}} , 
\end{equation} 
and for the solution, the following properties takes place
\begin{equation}\label{55}
 \|y^{n+1}\|_E = \|y^{n}\|_E, 
\quad E = I - \frac{1}{12} \tau^2 Q , 
 \quad n = 0,1,..., N-1 .   
\end{equation} 
\end{thm}

The condition (\ref{54}) ensures the positivity of the operator $E$.
To prove the equality (\ref{55}), it is sufficient to multiply equation (\ref{53}) by $\tau (y^{n + 1} + y^n)$.
Thus, we have the conservative property not in $H$ (see (\ref{15})), but only in $H_E$ (see (\ref{55})).
For the scheme (\ref{53}), we cannot prove the multiconservative property.

\section{Numerical experiments} 

The possibilities of the explicit and implicit schemes under the consideration for solving the advection equation are illustrated by the numerical results for the model problem.

\subsection{Model problem} 

We will consider the problem (\ref{1})-(\ref{3}) in the unit square:
\[
 \Omega = \{\bm x \ | \ \bm x= (x_1, x_2), \quad 0 < x_\alpha  <1, \quad \alpha  = 1,2\}.  
\] 
The initial condition is taken in the form
\[
 w^0(\bm x) = 2 \cdot 10^3 x_1^2 (1-x_1)^4  x_2^2 (1-x_2)^4 .
\] 
The components of the velocity $\bm v = (v_1, v_2)$ are defined by the stream function
\[
 \psi(\bm x) = \frac{1}{\pi} \sin(\pi x_1) \sin(\pi x_2) ,
\] 
so that
\[
 v_1 = \frac{\partial \psi}{\partial x_2},
 \quad  v_2 = - \frac{\partial \psi}{\partial x_1}.
\] 
The calculations are performed for $T = 5$.

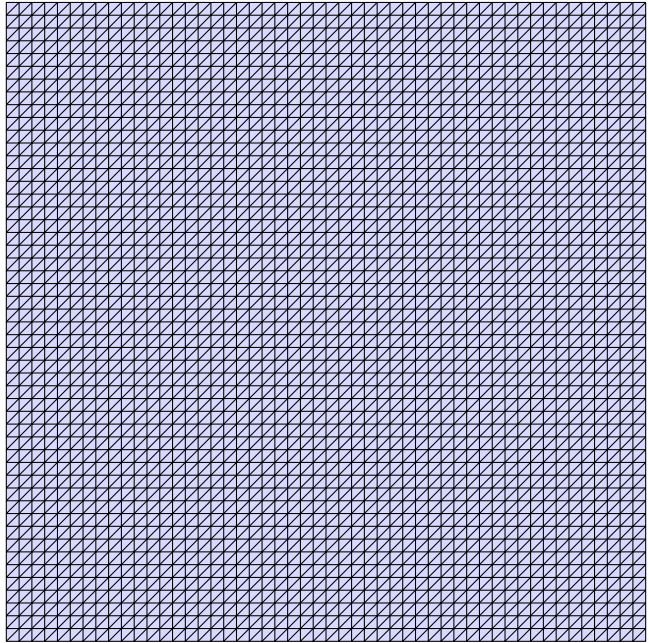
\begin{figure}[ht] 
  \begin{center}
    \begin{tikzpicture}[scale=0.85]
       \shade[top color=blue!15, bottom color=blue!15] (0,0) rectangle +(10,10);
	\draw[step=0.2] (0,0) grid (10,10);
	\foreach \x in {0,...,50} {
 	\draw (0.2*\x,0) -- (10,10-0.2*\x);
	}
	\foreach \y in {1,...,50} {
 	\draw (0,0.2*\y) -- (10-0.2*\y,10);
	}
    \end{tikzpicture}
    \caption{Grid with $h = 0.02$.} 
    \label{f-1}
  \end{center}
\end{figure}

The computational code  was implemented using the FEniCS \cite{fenics} numerical framework.
The finite element approximation in space is based on the use of continuous $P_1$ Lagrange element, namely, piecewise-linear elements.
A uniform grid is used for spatial domain.
The grid with the step $ h = 0.02 $ is shown in Fig.~\ref{f-1}.

\begin{figure}
\centering
\includegraphics[width=0.87\linewidth]{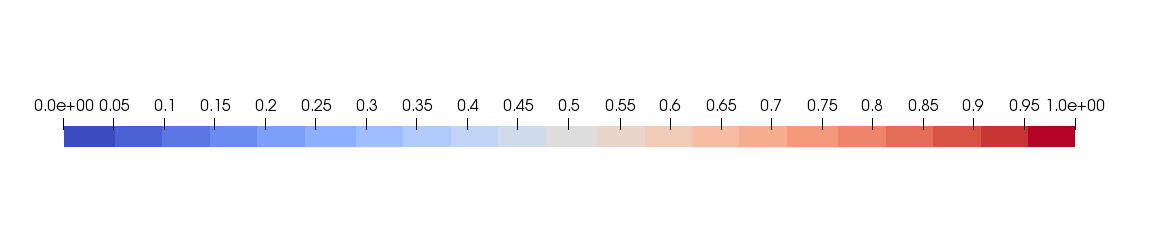} \\
\begin{minipage}{0.43\linewidth}
\centering
\includegraphics[width=\linewidth]{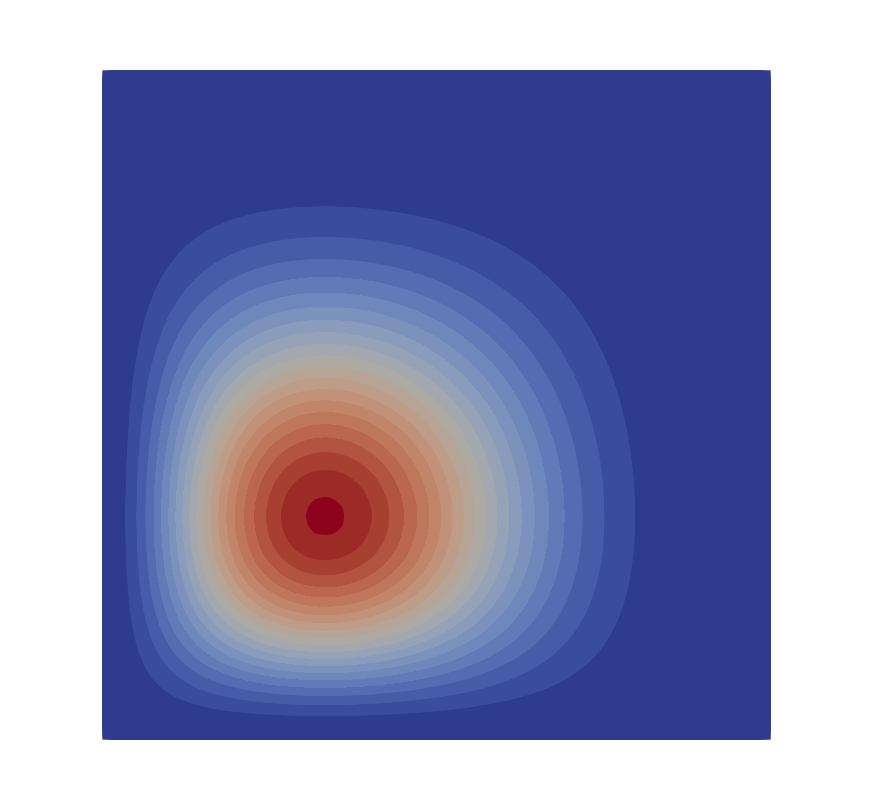}\\
$t=0$ \\
\includegraphics[width=\linewidth]{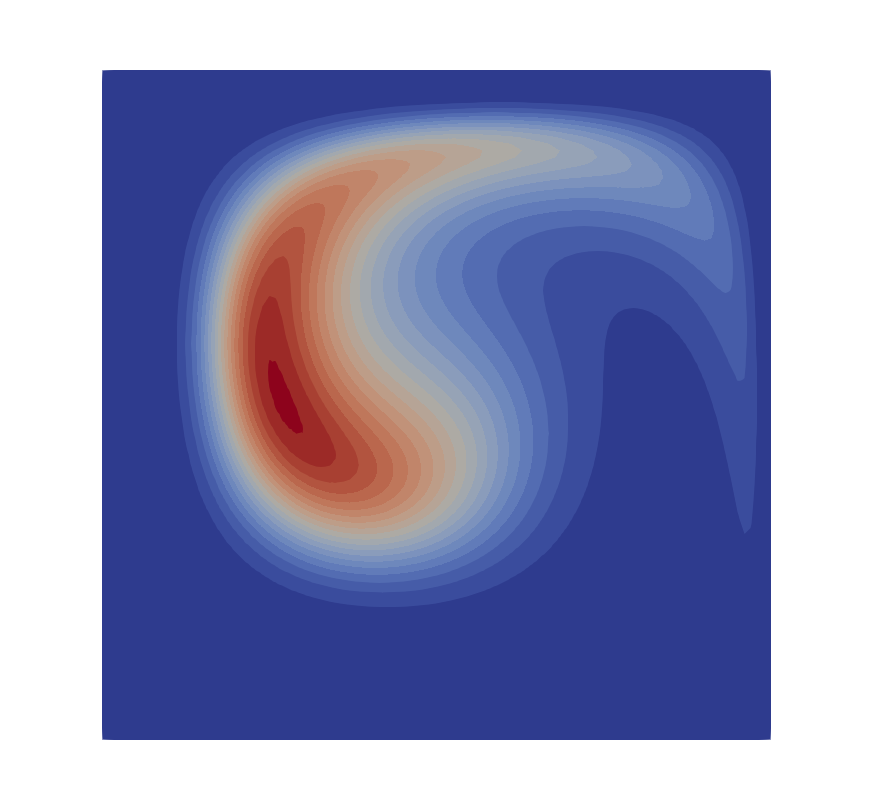}\\
$t=2$ \\
\includegraphics[width=\linewidth]{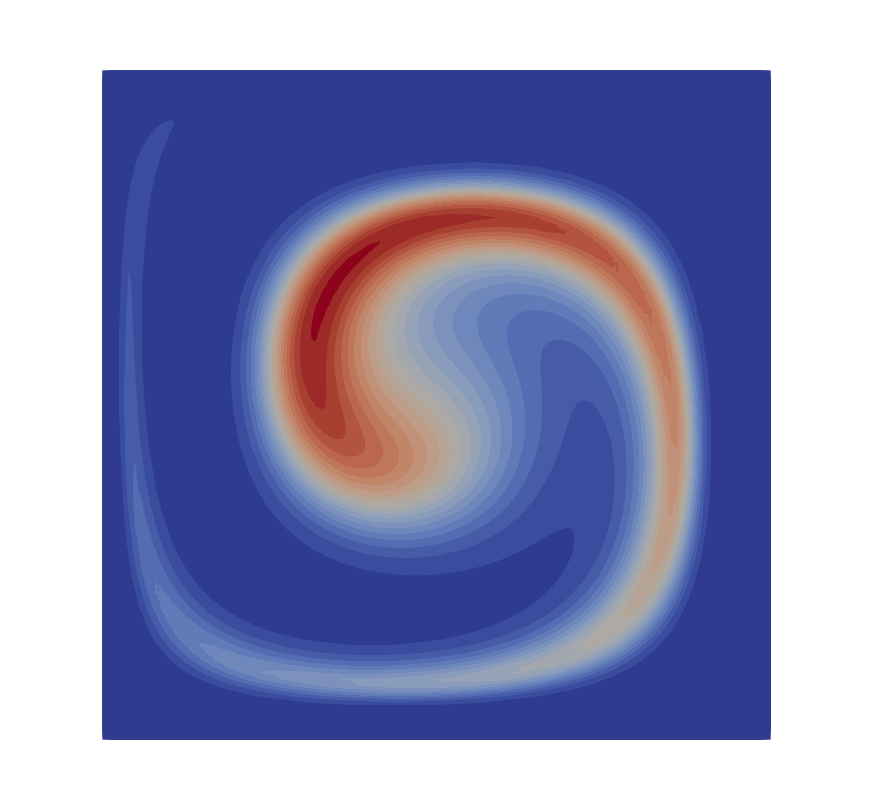}\\
$t=4$ 
\end{minipage}
\begin{minipage}{0.43\linewidth}
\centering
\includegraphics[width=\linewidth]{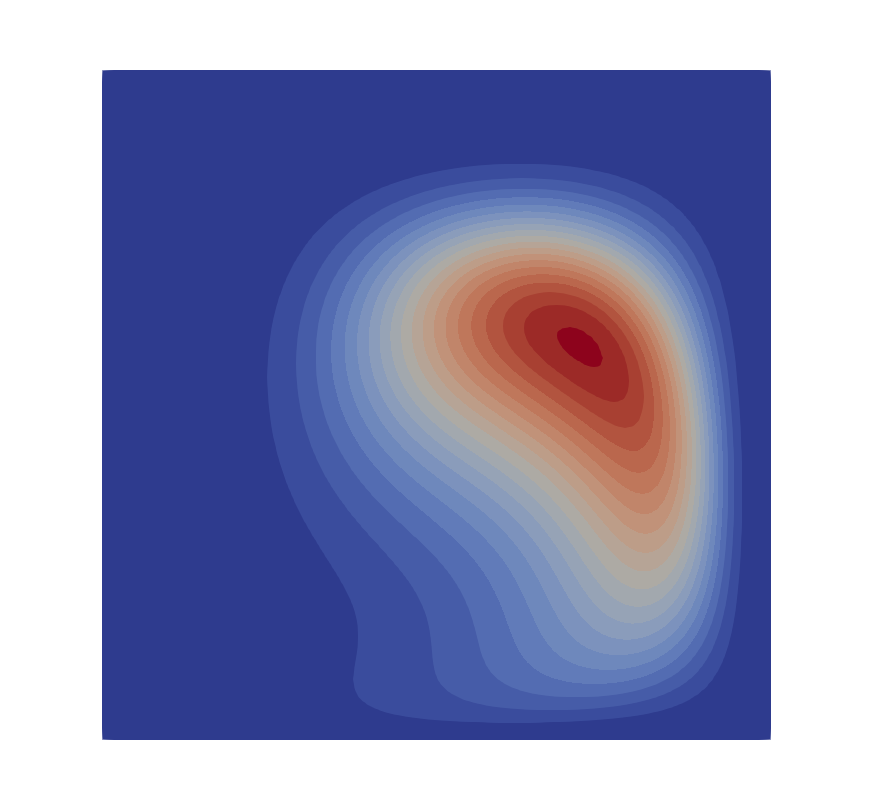}\\
$t=1$ \\
\includegraphics[width=\linewidth]{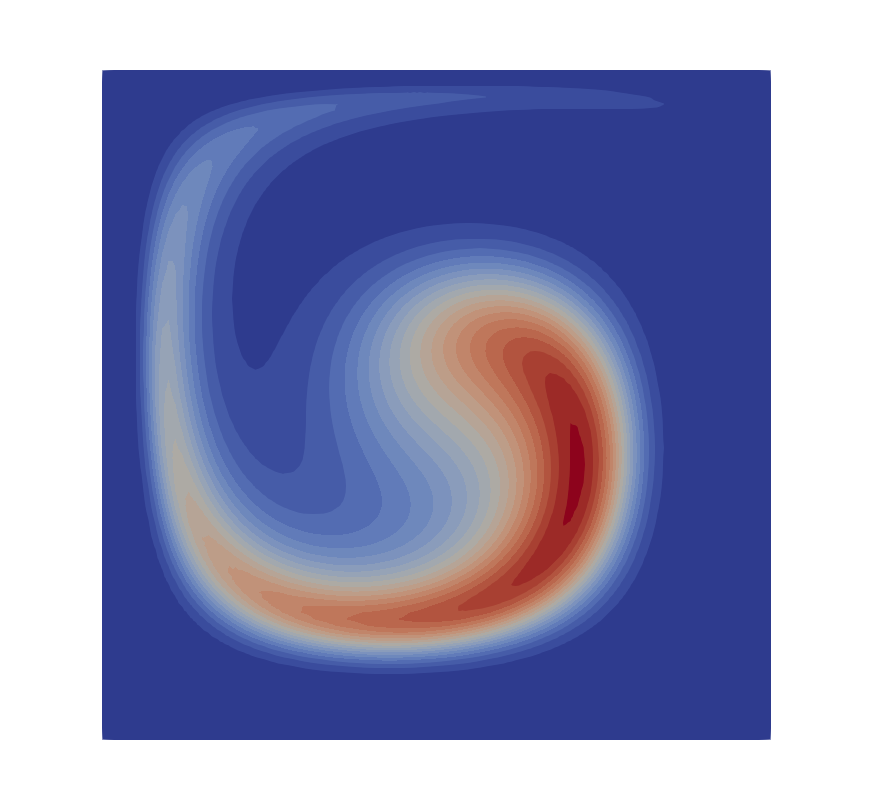}\\
$t=3$ \\
\includegraphics[width=\linewidth]{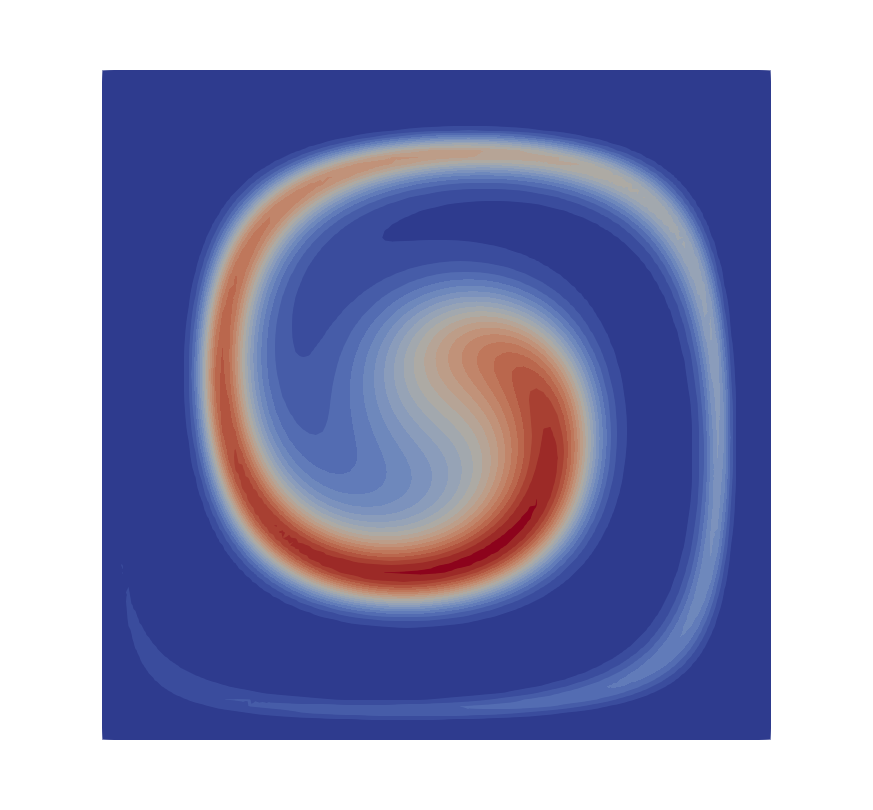}\\
$t=5$ 
\end{minipage}
\caption{The solution of the problem at different time-moments.}
\label{f-2}
\end{figure}

The accuracy of different approximations in time will be estimated by a reference solution.
It was obtained using the scheme under the consideration with an essentially small time step: $\tau = 2^{-11} \, 10^{-2}$.
The time-evolution of the solution on the grid with $h = 0.01$ (the main grid in space) is illustrated in Fig. \ref{f-2}.
For the relative error of the approximate solution, we have
\[
 \varepsilon (t) = \frac{\|y - \bar{y}\|}{\|\bar{y}\|} ,
\] 
where $\bar{y}$ is the reference solution.

\subsection{Explicit schemes} 

Prevously the explicit regularized schemes (\ref{18}), (\ref{30}) have been outlined.
For the choice of a regularizer in the form (\ref{34}), stability takes place with constraints (\ref{35}) on the time step.
It should be noted that we must orient on schemes with the operator (\ref{42}).

We present the results of calculations for the regularized scheme with $\beta = 1$.
The constraints for the step are related with the norm of the operator $A$.
Taking into account its skew-symmetry property, we have
\[
 \|A\| = \max_i |\lambda_i|,
\]     
where $\lambda_i, \ i =1,2, ..., N_h$ are the eigenvalues of the operator $A$:
\[
 A \varphi = \lambda \varphi .
\] 
Taking into account the representation (\ref{40}), this spectral problem corresponds to the spectral problem
\[
 K \psi = \lambda M \psi . 
\] 
Similarly, using the procedure of mass lumping with allowance for (\ref{42}), we get the spectral problem
\[
 K \psi = \lambda \widetilde{M} \psi . 
\] 

To solve the spectral problems with symmetrical matrices, we use the SLEPc library (Scalable Library for Eigenvalue Problem Computations) \cite{hernandez2005slepc}).
We apply the Krylov-Schur algorithm, a variation of the Arnoldi method, proposed by \cite{stewart2002krylov}.
The calculation results of the norm of the operator $A$ according to (\ref{40}) and (\ref{42}) on different grids are presented in the table~\ref{tab-1}.
It should be noted that $\| A \| = \mathcal{O}(h^{-1}) $ and the norm of the operator $A$ decreases approximately by two times when using mass lumping, because the maximum of permissible time step (see (\ref{35})) increases approximately twice.

\begin{table}[htp]
  \begin{center}
\caption{The norm of the operator $A$}
\label{tab-1}
  \begin{tabular}{lll}
   space grid ($h$) & $\|A\|$ for $M$ & $\|A\|$ for $\widetilde{M} $  \\
  \hline
   0.02  &  1.05288993e+02 & 5.59579462e+01  \\
   0.01  &  2.16001186e+02 & 1.14622718e+02  \\
   0.005 &  4.37491174e+02 & 2.31964151e+02  \\
  \end{tabular}
  \end{center}
\end{table} 

\begin{figure}
  \begin{center}
    \includegraphics[width=0.8\linewidth] {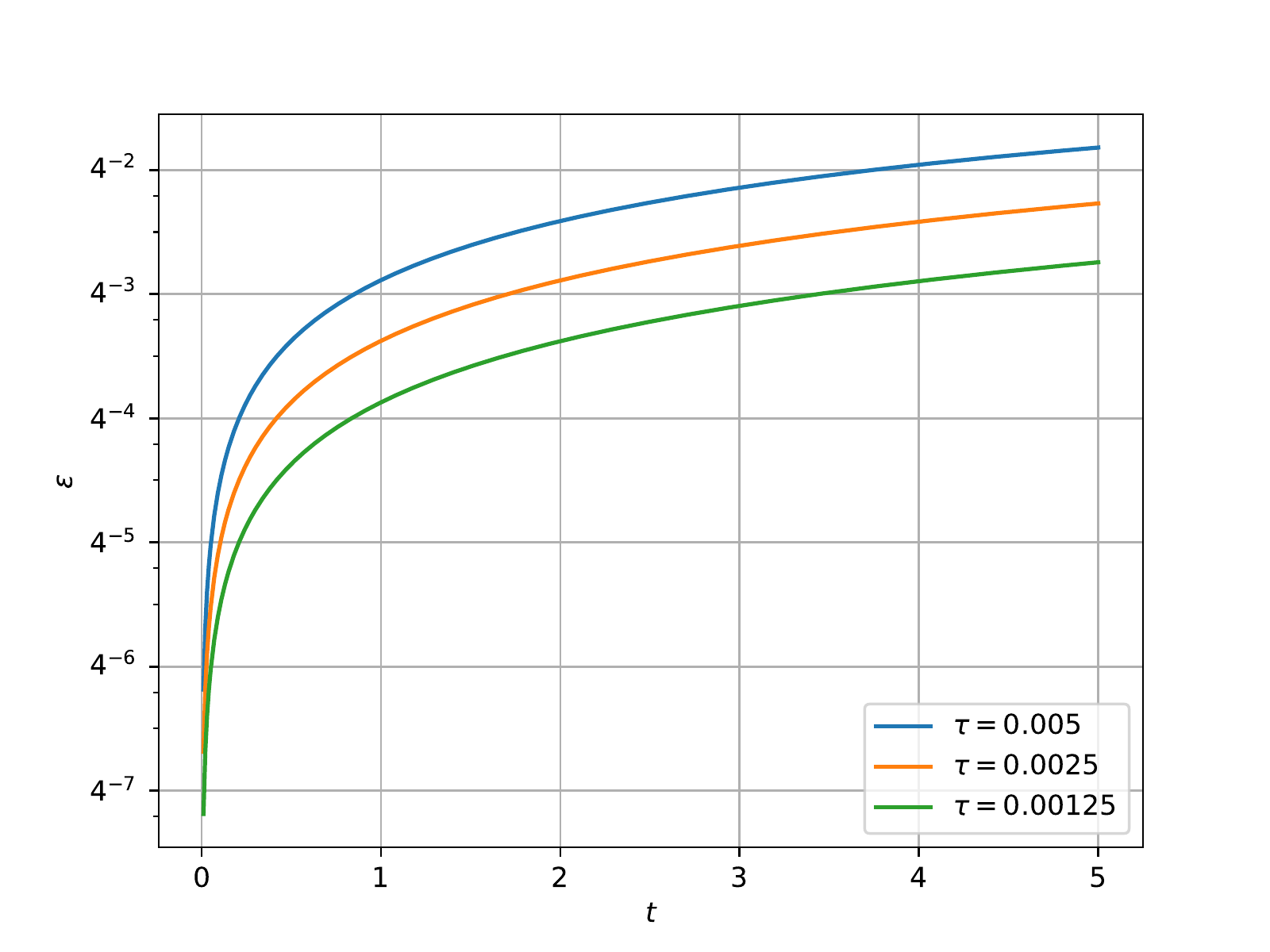}
	\caption{Error time-history of the regularized scheme.}
	\label{f-3}
  \end{center}
\end{figure}

The time-history of the error for various time steps for the regularized scheme (\ref{18}), (\ref{30}), (\ref{34}) with (\ref{42}) and $\beta = 1$ are shown in the Fig. \ref{f-3}. Convergence with first order in $\tau$ is observed.

Using the mass lumping procedure and taking into account (\ref{38}), the explicit Lax-Wendroff scheme is written as
\begin{equation}\label{56}
 \widetilde{M} \frac{z^{n+1} - z^n}{\tau} + K z^n + \frac{\tau}{2} G z^n= 0,
 \quad n = 0,1,..., N-1 . 
\end{equation} 
It corresponds to the regularized scheme (\ref{30}), whereh the operator $A$ is defined according to (\ref{42}), and for the operator $D$, we have
\begin{equation}\label{57}
 D = \widetilde{M}^{-1/2} G \widetilde{M}^{-1/2} .
\end{equation} 

The stability conditions of the scheme (\ref{30}) are formulated in theorem~\ref{t-4}.
The inequality $D \geq A^* A$ will be satisfied if the constant $\eta \geq 1$ in the inequality $D \geq \eta A^* A$.
We consider the spectral problem
\[
 A^*A \varphi = \lambda D \varphi , 
\] 
then $\eta = \lambda_{max}^{-1}$.
Taking into account (\ref{42}), (\ref{57}), this spectral problem corresponds the spectral problem
\begin{equation}\label{58}
 K^*\widetilde{M}^{-1} K \psi  = \lambda G \psi . 
\end{equation} 
The constraints on time step (\ref{33}) is related to the spectral problem (\ref{32}).
In the case (\ref{42}), (\ref{56}), it is equivalent to the spectral problem
\begin{equation}\label{59}
 G \widetilde{M}^{-1} G \psi  = \lambda (G - K^*\widetilde{M}^{-1} K ) \psi. 
\end{equation} 

The calculation results for the determination of $\eta$ and $\tau_0$ from the numerical solution of the spectral problems (\ref{58}), (\ref{59}) are given in table~\ref{tab-2}.
These data demonstrate the conditional stability of the Lax-Wendroff scheme, moreover $\tau_0 = \mathcal{O}(h)$.

\begin{table}[htp]
  \begin{center}
\caption{Parameters of the explicit Lax-Wendroff scheme}
\label{tab-2}
  \begin{tabular}{lll}
   space grid ($h$) & $\eta$ & $\tau_0$  \\
  \hline
   0.02  &  1.00098795e+00 & 1.73477111e-02  \\
   0.01  &  1.00025320e+00 & 8.47323207e-03  \\
   0.005 &  1.00006414e+00 & 4.17705891e-03  \\
  \end{tabular} 
  \end{center}
\end{table} 

\begin{figure}
  \begin{center}
    \includegraphics[width=0.8\linewidth] {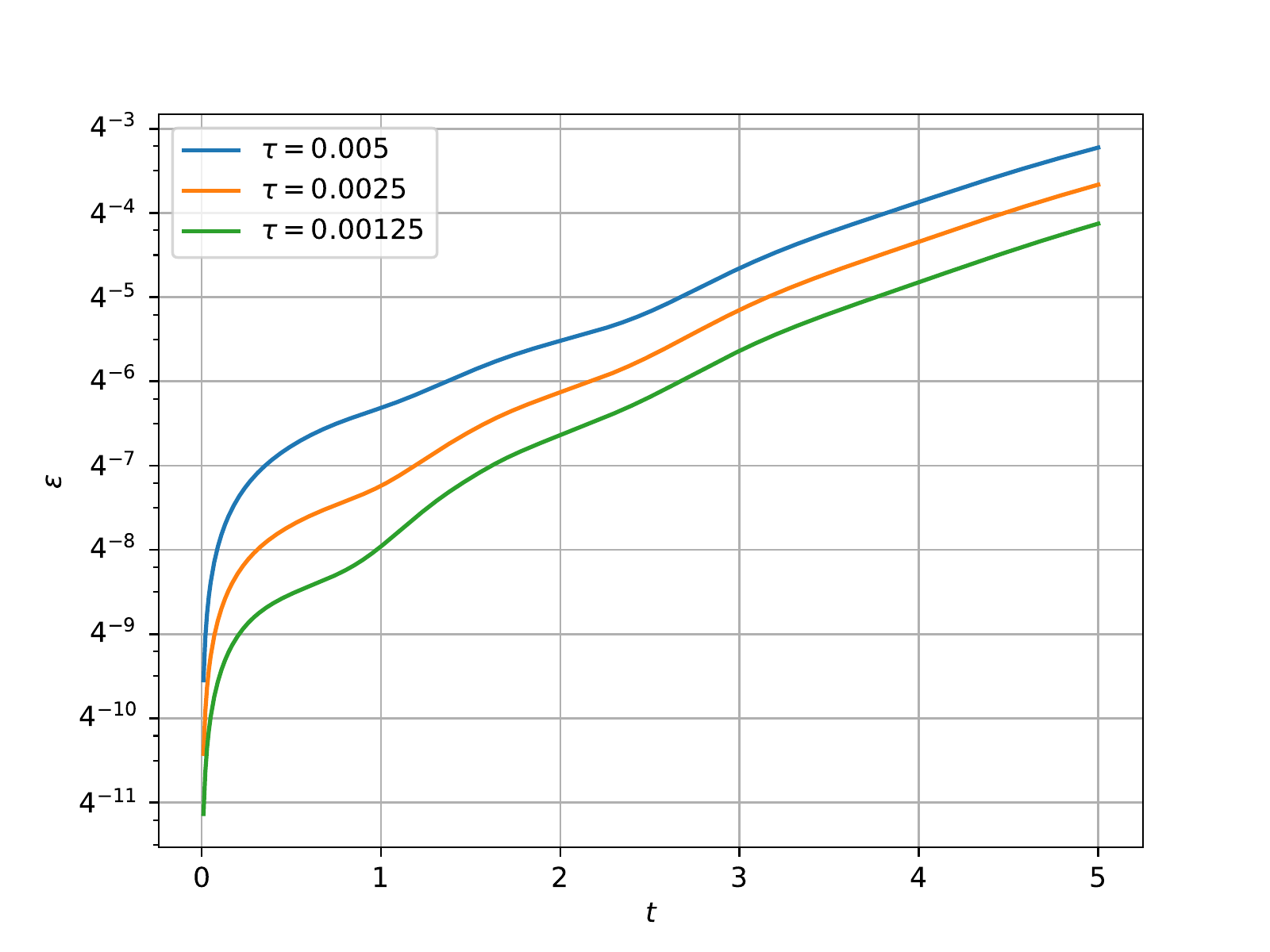}
	\caption{Error time-history of the Lax-Wendroff scheme.}
	\label{f-4}
  \end{center}
\end{figure}

We write the explicit Lax-Wendroff scheme (\ref{56}) in the form
\[
 \widetilde{M} \frac{z^{n+1} - z^n}{\tau} + K z^n - \frac{\tau}{2} K^*\widetilde{M}^{-1} K z^n +\frac{\tau}{2} R z^n  = 0,
 \quad n = 0,1,..., N-1 , 
\] 
Here the term $ R z ^ n $ with
\[
 R = G - K^*\widetilde{M}^{-1} K
\] 
distinguishes this scheme from the explicit Runge-Kutta scheme of the second-order accuracy in time.
In our problem, the convergence of the approximate solution with the second order in time is manifested (see Fig. \ref{f-4}), first of all, in the initial time interval, at which the smoothness of the solution is large enough and the influence of the term $R z^n$ is insignificant.
The first order of accuracy due to the term $R z^n$ begins to appear at $t \approx 1$.

\subsection{Implicit schemes}
 
\begin{figure}
  \begin{center}
    \includegraphics[width=0.8\linewidth] {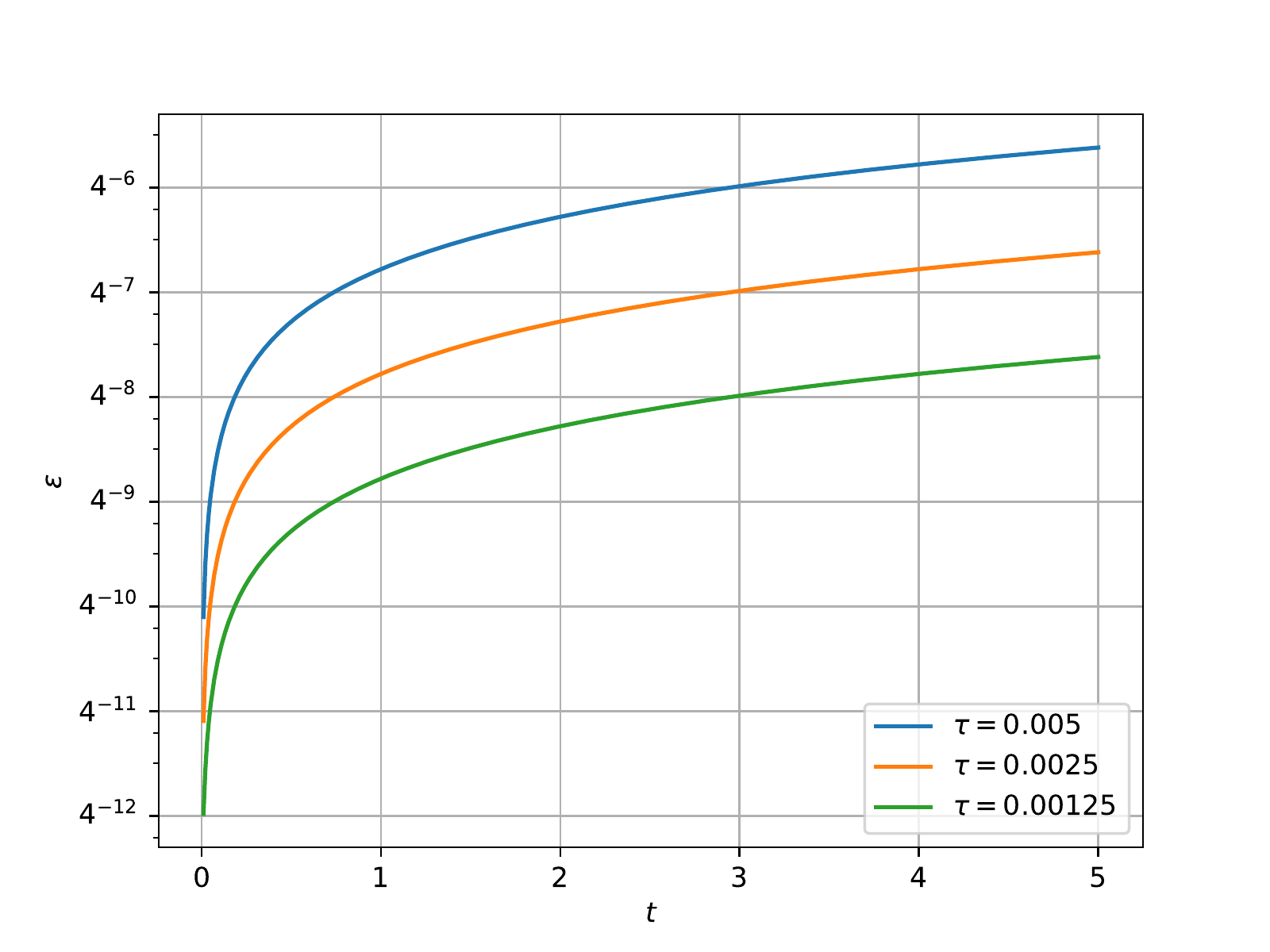}
	\caption{Error time-history of the Crank-Nicolson scheme.}
	\label{f-5}
  \end{center}
\end{figure}

In the class of implicit schemes, the Crank-Nicolson scheme (\ref{18}), (\ref{43}) seems to be the basic one.
Its accuracy in solving the model problem is illustrated in Fig. \ref{f-5}.
In comparison with the explicit regularized scheme (see Fig. \ref{f-3}), a higher accuracy of the approximate solution is observed, convergence is approximately of the second order in $\tau$.
In addition, this scheme stands out among all those considered schemes due to the fact that it is unconditionally stable.

\begin{figure}
  \begin{center}
    \includegraphics[width=0.8\linewidth] {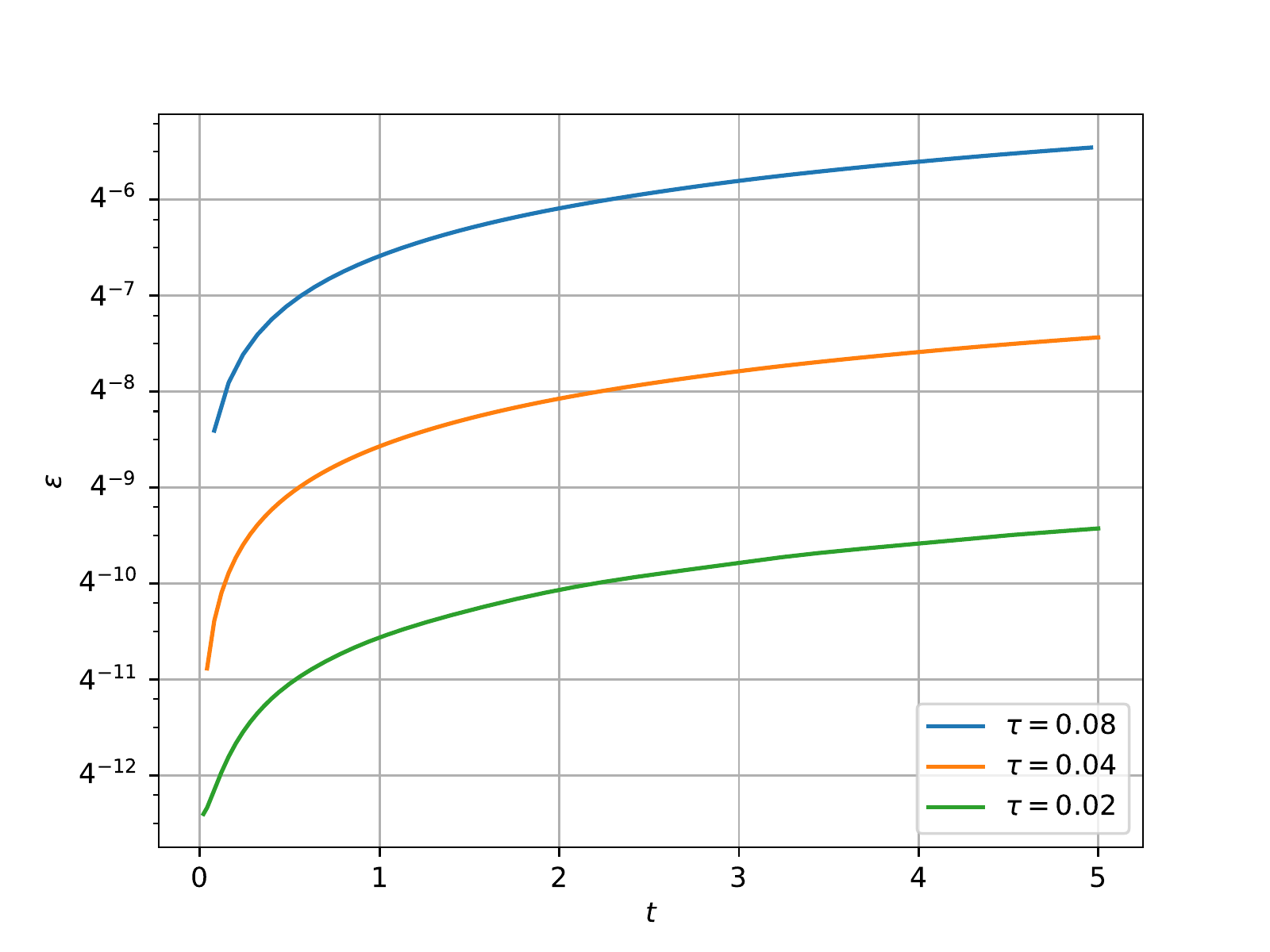}
	\caption{Error time-history of the scheme of higher accuracy.}
	\label{f-6}
  \end{center}
\end{figure}

Similar results for the scheme of higher order of accuracy are shown in Fig. \ref{f-6}.
The approximate solution using the lumping procedure is found from the equation (\ref{52}).
This scheme is absolutly stable and demonstrates the high accuracy of the approximate solution for substantially large time steps.

\begin{figure}
  \begin{center}
    \includegraphics[width=0.8\linewidth] {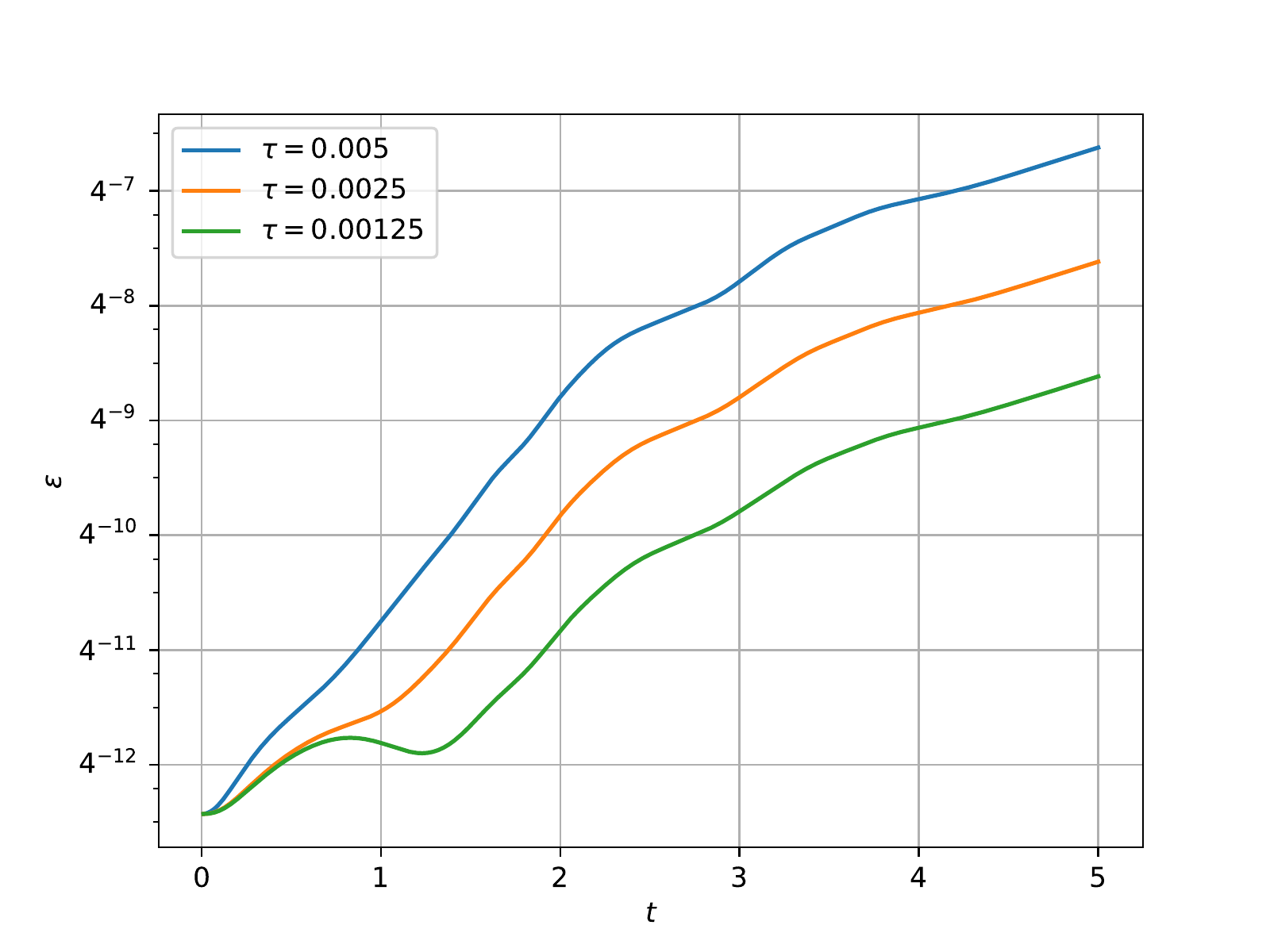}
	\caption{Accuracy of the implicit Lax-Wendroff scheme.}
	\label{f-7}
  \end{center}
\end{figure}

The accuracy of the implicit version of the Lax-Wendroff scheme is illustrated in Fig. \ref{f-7}.
As for the explicit scheme (see Fig. \ref{f-4}), there is a higher accuracy for the  initial time, when the solution is smooth in space.
Further, the effect of $Q-A^*A$ occures that results in the difference between the implicit Lax-Wendroff scheme (\ref{53}) and the fourth-order accuracy scheme (\ref{44}).
The scheme is stable for $\tau < \tau_0$, where $\tau_0$ corresponds to the right side (\ref{54}).
Numerical results for $\|Q\|$ and $\tau_0$, which are obtained from the solution of the spectral problem for the operator $\|Q\|$, are given in table \ref{tab-3} and show that $\|Q\| = \mathcal{O}(h^{-2})$ and $\tau_0 = \mathcal{O}(h)$.

\begin{table}[htp]
  \begin{center}
\caption{Parameters of the implicit Lax-Wendroff scheme}
\label{tab-3}
  \begin{tabular}{lll}
   space grid ($h$) & $\|Q\|$ & $\tau_0$  \\
  \hline
   0.02  &  3.22933843e+04 & 1.92767512e-02  \\
   0.01  &  1.33745164e+05 & 9.47221570e-03  \\
   0.005 &  5.44513748e+05 & 4.69446600e-03  \\
  \end{tabular} 
  \end{center}
\end{table} 

\section{Conclusions} 

\begin{enumerate}
\item In numerical simulation of advection processes, it is necessary to focus on writing the equation in a symmetric form (the half-sum of the advection operator in the divergent (conservative) form and the advection operator in the non-divergent (characteristic) form).
 In this case the advection operator is skew-symmetric.
 For the solution of the non-stationary advection equation, the multiconservative property holds, which is associated with the fulfillment of the set of conservation laws.
 \item The stability of known and new two-level difference schemes is investigated for the approximate solution of the Cauchy problem for the advection equation.
 Investigation of stability is carried out on the basis of the general theory of stability (well-posedness) of operator-difference schemes.
 The standard finite-element approximation is used in space, which ensures the skew-symmetry of the discrete advection operator.
 \item The standard explicit scheme for the advection equation belongs to the class of absolutely unstable ones.
 The explicit regularized scheme of the first order of accuracy is proposed, and the stability conditions are formulated.
 The stability conditions are established for the finite-element version of the classical explicit Lax-Wendroff scheme. For finite-element approximation in space, it is necessary to use various diagonalization algorithms of the mass matrix.
 \item For solving advection problems, the Crank-Nicolson scheme seems to be very attractive.
 It demonstrates the second order of accuracy and belongs to the class of unconditionally stable schemes.
 In addition, it demonstrates many conservation laws.
 The possibilities of using schemes with higher order of accuracy, which is also unconditionally stable and multiconservative, are also considered.
 The implicit version of the Lax-Wendroff scheme is proposed, and conditions for its stability are formulated.
 \item The properties of the explicit schemes are illustrated by the solution of the model two-dimensional problem for the advection equation. The focus is on studing the accuracy of various approximations in time.
\end{enumerate}

\section*{Acknowledgements}

This work is supported by the mega-grant of the Russian Federation Government (\#~14.Y26.31.0013).

\end{document}